\documentclass[12pt,a4paper,twoside]{article}
\usepackage{amsfonts} 
\usepackage{amssymb}
\usepackage{amsmath}

\usepackage{latexsym}
\usepackage{enumerate}
\usepackage[dvips]{graphicx}

\newcommand{\mode}{\mathrel{\text{mod}^*}}
\newtheorem{thm}{Theorem}[section]
\newtheorem{lem}{Lemma}[section]

\newtheorem{defi}{Definition}[section]

\newenvironment{dem}{\medbreak\noindent{\sc Proof: }}%
{\hfill $\diamond\diamond\diamond$\par\medbreak}

\def\mff{\mathfrak{f}}
\def\Ocal{\mathcal{O}}
\DeclareMathOperator{\Log}{Log}

\title{Improving on the Brun-Titchmarsh Theorem\footnote{AMS classification:
    11N13, 11N35, 11N36}\footnote{Keywords: Brun-Titchmarsh Theorem}}
\author{Oliver Ramar\'e \\ Jan Christoph Schlage-Puchta}
\date{2008}

\begin{document}

\begin{abstract}
  We prove that the number of primes in an interval of length $N$ is at
  most $2N/(\Log N+3.53)$ when $N$ is large enough. 
  This is obtained through a sieving process which can be seen as a hybrid between
  the large sieve and the Selberg sieve, and draws on what we call "local
  models".

  This paper has been published in
  \smallskip

 \noindent {\sc Acta Arithmetica in 2008, vol.
  131, nb. 4, pages 351--366}.
\end{abstract}

\maketitle

\section{Introduction}

The main result of this paper is the following Theorem:
\begin{thm}
  There exists an $N_0$ such that for all $N\ge N_0$ and all $M\ge1$
  we have
  \begin{equation*}
    \pi(M+N)-\pi(M)\le \frac{2N}{\Log N+3.53}.
  \end{equation*}
\end{thm}
Such theorems have been termed
``Brun-Titchmarsh'' Theorem by Linnik in
\cite{Linnik*61}.  Indeed, Titchmarsh proved such a
theorem for $q=1$, with a 
$\Log\Log(N/q)$ term instead of the~2 above, to establish the asymptotic for the number of
divisors of the $p+1$, $p$ ranging through the primes; he used the method
of Brun.  The constant~2 (with a $o(1)$) appeared for the first time in
\cite{Selberg*49}. In this work, Selberg
also shows that the constant $2+o(1)$ is optimal
in the above, if we are to stick to sieve methods in a fairly general context.
He expanded this theory, now known as the ``parity principle'',
in~\cite{Selberg*72}.

It is thus of interest to try to qualify the $o(1)$ --- in $2+o(1)$.  The first
upper bound of the shape $2N/(\Log N +c)$ with an unspecified but very
negative $c$ is due to van Lint \& Richert in \cite{van-Lint-Richert*65}
though \cite{Selberg*49} mentions without any proof such a result around equation
$(6)$ therein.
Bombieri gave in \cite{Bombieri*71} the value $c=-3$ and
Montgomery \& Vaughan the value $c=5/6$ in \cite{Montgomery-Vaughan*73}. The
section~22 of ``lectures on sieves'' \cite{Selberg*91} gives a
proof of $c=2.81$, a proof from which we shall take several elements.

Our method goes through a sieving process which can be seen
as a hybrid between the large sieve and the Selberg sieve, and draws
on what we call "local models". With no further input this would lead
to Selberg's results, with a tiny saving in the computations required.
But having at our disposal a machinery that allows such a high level of sieving
(the moduli we consider are as large as $N$), we can handle the
remainder in a more detailled way. One can see this process as a
weighted large sieve inequality, with weights adapted to the problem.
We note that we shall in between be confronted with the problem of
majorizing a step function by a polynomial, and that this problem
appears to be numerically more difficult than meet the eye.
And prior to that, we shall also be forced to use numerical means for
want of
a closed expression for our majorant. As a consequence
of these two points, we are not in a position to assert that our
result is best possible, even if we restrict our attention to our own
method!

As a matter of notations, we denote by $\sigma(d)$ the sum of the (positive)
divisors of $d$ while, for any $r\ge0$, we set
 $\eta_r(k)=\prod_{p|k}(p^{r+1}+1)/(p-1)$.

\section{Hilbertian inequalities}
\label{Hilbertian inequalities}
Let us start with a complex vector space $\mathcal{H}$ endowed with
a hermitian product $[ f|g]$, left linear and right
sesquilinear.

The easiest exposition goes through a formal definition:

\begin{defi}
  By an \emph{almost orthogonal system} in $\mathcal{H}$, we mean a
  collection of three sets of datas
  \begin{enumerate}
  \item a finite family $(\varphi^*_i)_{i\in I}$ of points\footnote{The reader
      may wonder why we chose to call the members of this family with a
      star. It is  to be consistent and to avoid confusion with
      notations that will appear later on.},
  \item  a finite family
    $(M_i)_{i\in I}$  of positive real numbers,
  \item a finite family 
    $(\omega_{i,j})_{i,j\in I}$ of complex numbers with
    $\omega_{j,i}=\overline{\omega_{i,j}}$,
  \end{enumerate}
  all of them given so that
  \begin{equation}\label{err}
    \forall(\xi_i)_i\in\mathbb{C}^I,\quad
    \|\sum_i\xi_i\varphi^*_i\|^2
    \le
    \sum_iM_i|\xi_i|^2+\sum_{i,j}\xi_i\overline{\xi_j}\omega_{i,j}.
  \end{equation}
\end{defi}
We are to comment on this definition.
If the family $(\varphi^*_i)_{i\in I}$ were orthogonal, we could ask for
equality with $M_i=\|\varphi_i^*\|^2$. As it turns out, in applications we
have in mind, this family is not orthogonal, but almost so.
It is this almost orthogonality that the above condition is meant to measure.

Our first lemma reads as follows
\begin{lem}\label{lemPS1}
  For any finite family $(\varphi^*_i)_{i\in I}$ of points of $\mathcal{H}$,
  the system build with $M_i=\sum_j|[\varphi^*_i|\varphi^*_j]|$ and
  $\omega_{i,j}=0$ is almost orthogonal.
\end{lem}
So that, when $[\varphi^*_i|\varphi^*_j]$ is small for $i\neq j$, then $M_i$ is
indeed close to $\|\varphi^*_i\|^2$
\begin{dem}
  We write
  \begin{equation*}
    \Bigl\|\sum_i\xi_i\varphi^*_i\Bigr\|^2=
    \sum_{i,j}\xi_i\overline{\xi_j}[\varphi_i^*|\varphi_j^*]
  \end{equation*}
  and simply apply $2|\xi_i\overline{\xi_j}|\le |\xi_i|^2 +|\xi_j|^2$.
  The lemma readily follows.
\end{dem}

Here is an enlightening reading of this lemma: the hermitian form that
appears has a matrix whose diagonal terms are the $\|\varphi^*_i\|^2$'s. A
theorem of Gershgorin says that all eigenvalues of this matrix are to lie in
the so-called \emph{Gershgorin's disc} centered on one $\|\varphi^*_i\|^2$ and
with radius $\sum_{j\neq i}|[\varphi^*_i|\varphi^*_j]|$.  This approach is due
to \cite{Elliott*71}.  It has a drawback: we do not know that each Gershgorin
disc does indeed contain an eigenvalue, a flaw that is somehow repaired in the
above lemma.

In general, and only under~\eqref{err}, we get the following kind of Parseval
inequality:  
\begin{lem}\label{lemPS} For any almost orthogonal system,
  and any $f\in\mathcal{H}$, let us set $\xi_i=[ f|\varphi^*_i]/M_i$. We
  have 
  \begin{equation*}
    \sum_i M^{-1}_i|[ f|\varphi^*_i]|^2\le \|f\|^2
    +\sum_{i,j}\xi_i\overline{\xi_j}\omega_{i,j}.
  \end{equation*}
\end{lem}
Once again, the orthogonal case is enlightening: when the $(\varphi^*_i)$ are
orthogonal, then we may take $M_i=\|\varphi^*_i\|^2$ and
$\omega_{i,j}=0$. The LHS becomes the 
square of the norm of the orthonormal projection of $f$ on the subspace
generated by the $\varphi^*_i$'s. 

Without the $\omega_{i,j}$'s and appealing to Lemma~\ref{lemPS1}, this is due to
Selberg, as mentioned in section~2 of \cite{Bombieri*87}
and in \cite{Bombieri*71}.
\begin{dem}
  For the proof, simply write
  \begin{equation*}
    \bigl\|f-\sum_i\xi_i\varphi^*_i\bigr\|^2\ge0
  \end{equation*}
  and expand the square. We take care of 
  $\|\sum_i\xi_i\varphi^*_i\|^2$ by using \eqref{err}, getting
  \begin{equation*}
    \|f\|^2-2\Re \sum_i\overline{\xi_i}[ f|\varphi^*_i]
    +\sum_iM_i|\xi_i|^2+\sum_{i,j}\xi_i\overline{\xi_j}\omega_{i,j}\ge0.
  \end{equation*}
  We now choose $\xi_i$'s to the best of our interest, neglecting the bilinear form
  containing the $\omega_{i,j}$'s. We take
  $\xi_i=[ f|\varphi^*_i]/M_i$. The lemma readily follows.
\end{dem}
Combining Lemma~\ref{lemPS} together with Lemma~\ref{lemPS1} yields
what is usually known as ``Selberg's lemma'' in this context. The
introduction of the $\omega_{i,j}$'s is due to the authors to enable a
refined treatment of the error term as well as an hybrid way between
weighted large sieve results and Selberg sieve results.

The value of $\xi_i$ in the statement is usually of no importance, only its
order of magnitude being relevant.

We now discuss a special problem to introduce our next Theorem.
In some cases, a partial treatment of the bilinear form is readily available in
the shape of
\begin{equation}\label{erreurcond}
  \forall(\xi_i)_i\in\mathbb{C}^I,\quad
  \Bigl\|\sum_i\xi_i\varphi^*_i\Bigr\|^2
  \le
  \sum_i M_i|\xi_i|^2+\bigg(\sum_{i}|\xi_i|n_i\bigg)^2
  +\sum_{i,j}\xi_i\overline{\xi_j}\omega_{i,j}
\end{equation}
for some positive $M_i$, and $n_i$ (here again, $M_i$ is generally an
approximation of $\|\varphi^*_i\|^2$). 
With such an inequality at hand, the above proof leads to
\begin{equation}\label{source}
    \|f\|^2-2\Re \sum_i\overline{\xi_i}[ f|\varphi^*_i]
    +\sum_iM_i|\xi_i|^2+\bigg(\sum_{i}|\xi_i|n_i\bigg)^2
    +\sum_{i,j}\xi_i\overline{\xi_j}\omega_{i,j}\ge0.
\end{equation}
When using it, we take for
$\varphi^*_i$ a kind of local approximation of $f$, which implies that
we can assume $[
f|\varphi^*_i]$ to be a non-negative real number. 
It is then readily seen that the $\xi_i$'s that minimize the RHS are non-negative.
Finally, we are led to choose these $\xi_i$'s so as to minimize
\begin{equation*}
    \|f\|^2-2\sum_i \xi_i[ f|\varphi^*_i]
    +\sum_iM_i\xi_i^2+\bigg(\sum_{i}\xi_i n_i\bigg)^2.
\end{equation*}

We handle the optimization of~\eqref{source} with calculus
by setting $\xi_i=\zeta_i^2$. After some manipulations, we conclude that there
exists a subset $I'$ of $I$ such that $\xi_i=0$ if
$i\in I\setminus I'$ and
\begin{equation}
  \label{eq:22}
  \forall i\in I',\quad\xi_i=\frac{[ f|\varphi^*_i]-Xn_i}{M_i},
  \quad
  X=\frac{\sum_{j\in I'}n_j[ f|\varphi^*_j]/m_j}
  {1+\sum_{j\in I'}n_j^2/m_j}
\end{equation}
provided
\begin{equation}
  \label{eq:2}
  \forall i\in I',\quad[ f|\varphi^*_i]/n_i\ge X.
\end{equation}
However, determining optimal $I'$ is difficult: the index $i$ appears
on the left-hand side of~\eqref{eq:2}, but also on its right-hand side
since the definition of $X$ depends on whether this index
belongs to $I'$ or not.
It is easier to set
\begin{equation}
  \label{eq:51}
  \xi_i=\frac{[ f|\varphi^*_i]-Yn_i}{M_i},
 \end{equation}
for a $Y$ to be chosen but which guarantees $\xi_i\ge0$. The optimal
$Y$ is of course $Y=X$. 

Once we have inferred the form of these weights, we can simply plug
them in the proof of Lemma~\ref{lemPS} without even mentioning~\eqref{erreurcond}.
Here is the theorem we have reached:

\begin{thm}\label{thmweighted}
  Let an almost orthogonal system be given with notations as above
  and let $f\in\mathcal{H}$. Let also $Y$ be a non-negative real number
  and $(n_i)_i$ be non-negative real numbers.
  Assume that
  $[ f|\varphi^*_i]$'s
  are real numbers. Set
  $\xi_i=([ f|\varphi^*_i]-Yn_i)/M_i$  for all~$i$. Then we have
  \begin{equation*}
    \sum_i 
    M_i\xi_i^2+2Y\sum_i n_i\xi_i
    -\sum_{i,j}\xi_i\overline{\xi_j}\omega_{i,j}\le \|f\|^2.
  \end{equation*}
\end{thm}

Of course, the preliminary discussion tells us that it will be better to have
$\xi_i\ge0$, but the statement is valid as is, and may offer some more flexibility.

\section{Integers coprime to a fixed modulus in an interval}
\label{sec3}

Let $\mff$ be a positive integer and let us define by $\rho=\phi(\mff)/\mff$.
We study here the following two
functions of the real non-negative variable $u$: 
\begin{equation*}
        \left\{
        \begin{array}{rcl}
        \theta_{\mathfrak{f}}^-(u)&=&    \displaystyle 
    \min_{y\in\mathbb{R}}
       \min_{\substack{0\le x\le u\\ x\in\mathbb{R}}} \biggl(
      \sum_{\substack{y<n\le y+x,\\ (n,\mathfrak{f})=1}}1
      -\rho x\biggr)
    ,\\ \displaystyle
    \theta_{\mathfrak{f}}^+(u)&=&\displaystyle
    \max_{y\in\mathbb{R}}
       \max_{\substack{0\le x\le u\\ x\in\mathbb{R}}} \biggl(
      \sum_{\substack{y<n\le y+x,\\ (n,\mathfrak{f})=1}}1
      -\rho x\biggr).
        \end{array}
        \right.
  \end{equation*}

The introduction of these two functions is inspired from section~22
of ``lectures on sieves'' in \cite{Selberg*91}. 
In order to compute them, we need to restrict both $x$ and $y$ to
integer values. This is the role of next lemma.

  \begin{lem}\label{remainder}
    We have
    \begin{equation*}
        \left\{
        \begin{array}{rcl}
        \theta_{\mathfrak{f}}^-(u)&=&    \displaystyle 
    \min_{\ell\in\mathbb{N}}\left(
       \min_{\substack{k\in\mathbb{N},\\ 0\le k\le u}}
       \biggl(
       \sum_{\substack{\ell+1 \le n\le \ell+k-1,\\ (n,\mathfrak{f})=1}}1
      -\rho k\biggr),
      \sum_{\substack{\ell+1\le n\le \ell+[u],\\ (n,\mathfrak{f})=1}}1-\rho u
      \right)
    ,\\ \displaystyle
    \theta_{\mathfrak{f}}^+(u)&=&\displaystyle
    \max_{\substack{k,\ell\in\mathbb{N},\\ k<u+1}}
       \biggl(
      \sum_{\substack{\ell\le n\le \ell+k-1,\\ (n,\mathfrak{f})=1}}1
      -\rho (k-1)\biggr)
        \end{array}
        \right.
  \end{equation*}
  The function $\theta_{\mathfrak{f}}^+$ is a non-decreasing step function which is left
  continuous with jump points at integer points. The function
  $\theta_{\mathfrak{f}}^-$ is non-increasing continuous : it alternates
  linear pieces of directing coefficient $-\rho$ and constant pieces. Changes
  occur at integer points. Both are constant if $u\ge\mff$.
\end{lem}

\begin{dem}
  We start with $\theta^+_\mff$. First fix $y$. The function $\sum_{y<n\le
  y+x}w(n)-\rho x$ is linear non-increasing in $x$ from 0 to $1-\{y\}$, then from
  $1-\{y\}$ to $2-\{y\}$ and so on. Its maximum value is reached at $x=0$ or
  $x=k-\{y\}$ for some integer $k$, thus
  \begin{equation*}
    \theta_{\mathfrak{f}}^+(u)
    =\max_{y\in\mathbb{R}}
    \max_{\substack{k\in\mathbb{N},\\ k\le u+\{y\}}}
    \biggl(
    \sum_{\substack{y<n\le [y]+k}}w(n)+\rho(-k+\{y\})
    \biggr).
  \end{equation*}
  The condition $k\le u+\{y\}$ is increasing in $\{y\}$ and so is the term to maximize. We
  may take thus $y$ to be just below an integer $\ell$, reaching the expression we
  announced.

  As for $\theta^-_\mff$, we start similarly by fixing
  $y$. Minimum is reached at $x=k-\{y\}-0$ or at $x=u$, where $k$ is an integer and the
  $-0$ means we are to take $x$ just below this value. We get
  that $\theta_{\mathfrak{f}}^-(u)$ equals
  \begin{equation*}
    \min_{y\in\mathbb{R}}\Biggl(
      \min_{\substack{k\in\mathbb{N},\\ k\le u+\{u\}}}
      \biggl(\sum_{y< n\le [y]+k-1}w(n)+\rho(\{y\}-k)\biggr),
      \sum_{y< n\le [y]+u}w(n)-\rho u
    \Biggr).
  \end{equation*}
  As far as the last sum is concerned, the worst case is when $y$ is an
  integer $\ell\ge0$, so it reduces to
  \begin{equation}\label{already}
    \min_{\ell\in\mathbb{N}}\biggl(
    \sum_{\ell+1\le n\le\ell+u}w(n)-\rho u\biggr).
  \end{equation}
  For the first minimum, we distinguish between $k\le [u]$ and $k=[u]+1$
  (which can only happen if $u$ is \emph{not} an integer). If $k\le [u]$, we
  may take $y$ to be integral. If $k=[u]+1$, then $\{y\}\ge 1-\{u\}$ which is
  indeed the worst case: we take $y=\ell+1-\{u\}$. This last contribution
  turns out to be exactly the same as the one in \eqref{already}. 
\end{dem}

Next we consider the function
\begin{equation}
  \theta^*_\mff(v)=\max(\theta^+_\mff(1/v),-\theta^-_\mff(1/v))
\end{equation}
which  is right continuous with jump points at the $1/m$'s, where $m$ 
ranges over the integers from $1$ to $\mff$. Of course, $\theta^*_\mff(1)=\theta^+_\mff(1)=1$.

\subsubsection*{Case of $\mff=210$}














Here is our function:
\begin{equation*}
  \theta^*_{210}(1/u)=
  \begin{cases}
    1&\text{if $0< u\le 1$}\\
    54/35&\text{if $1< u\le 3$}\\
    57/35&\text{if $3< u\le 7$}\\
    76/35&\text{if $7< u\le 9$}\\
    79/35&\text{if $9< u\le 79/8$}\\
  \end{cases}
  \begin{cases}
    8u/35&\text{if $79/8\le u\le 10$}\\
    16/7&\text{if $10< u\le 13$}\\
    82/35&\text{if $13< u\le 17$}\\
    94/35&\text{if $17< u\le 41/2$}\\
    8u/35-2&\text{if $41/2\le u\le 22$}\\
    106/35&\text{if $22< u\le 210$}
  \end{cases}
\end{equation*}

The following plot displays the step function $\theta^*_{210}$ as well as the
optimizing polynomial we shall compute in section~\ref{otp}.

\begin{figure}[!htb]
\centering\includegraphics[bb=0 50 576 771, clip, scale=0.48, angle= -90]{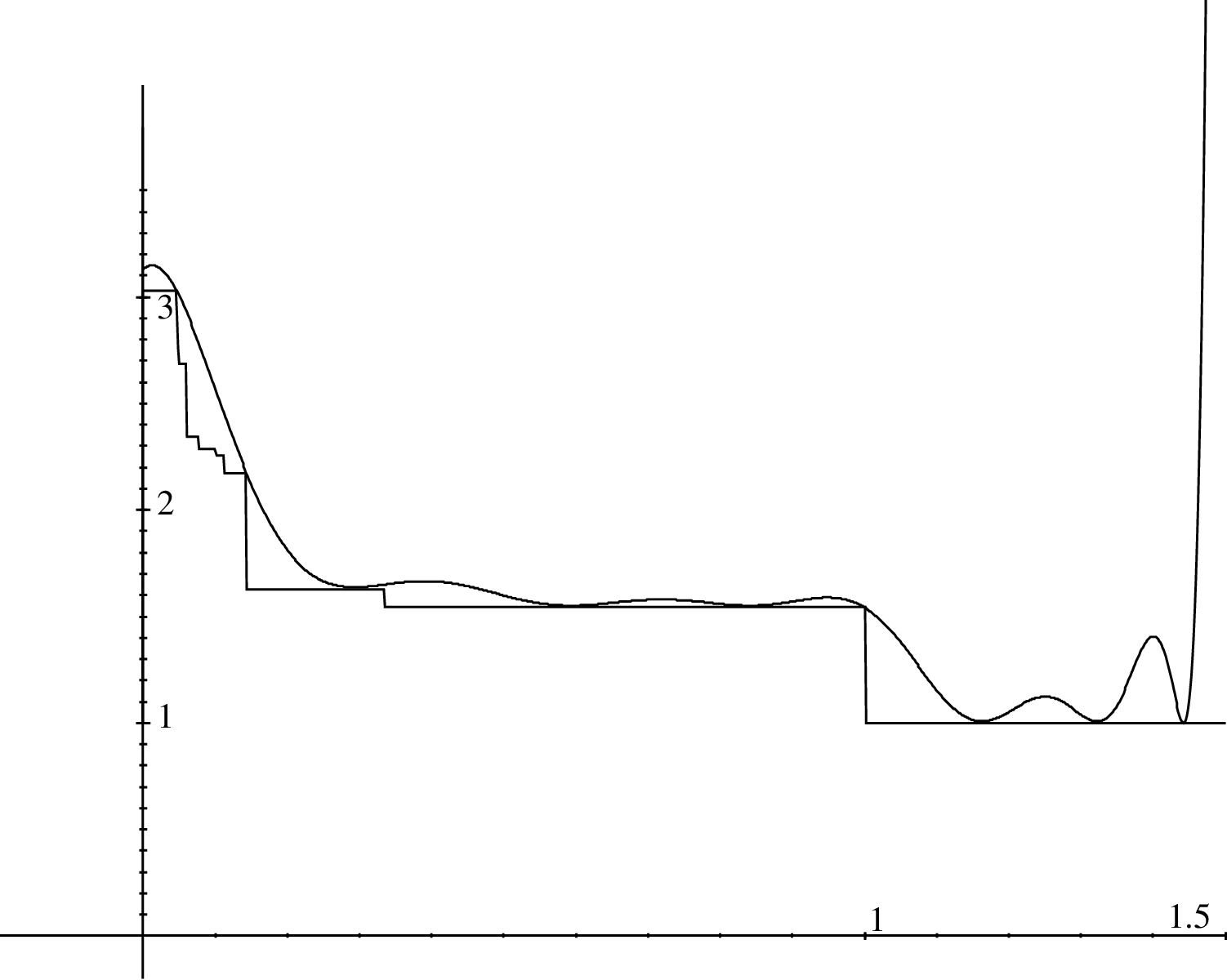}
\caption{Comparison of $\theta^*_{210}$ and the optimizing polynomial}
\end{figure}


\subsubsection*{Polynomial approximation of $\theta^*_\mff(v)$}

We shall  require a good polynomial upper bound for $\theta*_\mff$:
\begin{equation}\label{upb}
  \theta^*_\mff(v)\le \sum_{0\le r\le R} b_r v^r.
\end{equation}
Finding such an approximation turned out to be much more tricky than expected.
Our first idea has been to start with a polynomial approximation of $\theta^*_\mff(v)$ of the form
\begin{equation*}
  \bigr|\theta^*_\mff(v)-\sum_{0\le r\le R}\tilde{b}_r v^r\bigr|\le \epsilon
\end{equation*}
for $0\le v\le V$ from which an upper bound is easily derived by
increasing the constant term.
We carried out this scheme with Bernstein polynomials, with poor
numerical results and fitting. We even tried to achieve such an
approximation on a larger interval since endpoints are notoriously
troublesome: this helped a bit but not by much, despite the fact that we
used polynomials of very high degree (up to 200). We finally decided
for a different scheme exposed in the last section.


\section{Local models for the sequence of primes}

\subsection{Choice of the local system}
\label{Choice of the local system}
Let us start with a general discussion on what ``sieving'' means. Sieving is about
gaining information on a sequence from what we know of it modulo~$d$ for
several~$d$'s. If one looks at the sequence of primes modulo~$d$ and if we
neglect the prime divisors of $d$, it simply is the set of reduced residue
classes modulo~$d$. Thus, on one hand we have the
characteristic function of the primes of the interval $[M+1,M+N]$, say $f$, and on
the other hand the characteristic function $\varphi_d$ of the integers in this
interval that are coprime to $d$ for all $d\le\sqrt{N}$. Notice here that
it is enough to restrict our attention to squarefree $d$'s.

On recalling what we did in section~\ref{Hilbertian inequalities}, we could
simply try to get an approximation of $f$ in terms of the $\varphi_d$'s.
However, the study there is patterned for almost orthogonal $\varphi_q$'s,
which is not the case of the sequence $(\varphi_d)_d$: if $q|d$, knowing that
a given integer is coprime with $d$ implies it is coprime with~$q$, so there
is redundancy of information. It implies in turn that these functions are far
from being linearly independent. We unscrew the situation in the following
way. When $d$ is squarefree, we set
\begin{equation}
  \label{eq:10}
  \frac{d}{\phi(d)}\varphi_d=\sum_{q|d}\varphi_q^*
\end{equation}
where
\begin{equation}
  \label{eq:12}
  \varphi_q^*(n)=\mu(q)c_q(n)/\phi(q)
\end{equation}
and $c_q(n)$ is the Ramanujan sum given by
\begin{equation}
  \label{eq:13}
  c_q(n)=\sum_{a\mode q}e(na/q)=\sum_{\ell|q}\ell\mu(q/\ell).
\end{equation}
Verifying \eqref{eq:10} is easy:
\begin{equation*}
  \sum_{q|d}\mu(q)c_q(n)/\phi(q)=
  \prod_{p|d}\left(
    1-\begin{cases}
      1&\text{if $p|n$}\\
      -1/(p-1)&\text{otherwise}
    \end{cases}
  \right).
\end{equation*}

In our problem, we shall select a fixed integer $\mff$ that will be taken to
be $210$ at the end of the proof and consider the characteristic function $w$
of the points in $[M+1,M+N]$ that are coprime with $\mff$. This being chosen,
our hermitian product on sequences over $[M+1,M+N]$ is defined by
\begin{equation}
  \label{eq:15}
  [g|h]=\sum_{M+1\le n\le M+N}w(n)g(n)\overline{h(n)}.
\end{equation}
Furthermore, we take the moduli $q$ in the set
\begin{equation}
  \label{eq:16}
  \bigl\{q\ /\  \sigma(q)\le S,\mu^2(q)=1,(q,\mff)=1\bigr\},
\end{equation}
where $\sigma(q)=\sum_{d|q}d$. The reason for this choice will become apparent
later on. 

\subsection{Study of the local models}
\label{Study of the local models}

Notice that
\begin{equation}
  \label{eq:17i}
  [\varphi_q^*|\varphi_{q'}^*]
  =\frac{\mu(q)}{\phi(q)}\frac{\mu(q')}{\phi(q')}
  \sum_nw(n)c_q(n)c_{q'}(n).
\end{equation}
We note that when $q$ and $q'$ have a common factor, say $\delta$, then
$c_\delta(n)^2=\phi((n,\delta))^2$ would factor out: this contribution
is non-negative and we want to use this fact here. Let $\Delta$ be a
squarefree integer coprime with $\mff$.  Write $(q,q',\Delta)=\delta$,
so that $ [\varphi_q^*|\varphi_{q'}^*]$ equals
\begin{equation*}
  \frac{\mu(q)\mu(q')}{\phi(q)\phi(q')}
  \sum_{\substack{\ell|q/\delta\\\ell'|q'/\delta\\ h|\delta}}
  \ell\mu(q/\ell)\ell'\mu(q/\ell')(\mu\star\phi^2)(h)
  \left(\frac{\phi(\mff)}{\mff}\frac{N}{h[\ell,\ell']}
    +R_{h[\ell,\ell']}(M,N,\mff)\right)
\end{equation*}
where
\begin{equation}
  \label{eq:19}
  R_{d}(M,N,\mff)=
  \sum_{\substack{M+1\le n\le M+N\\ d|n}}w(n)
  -\frac{\phi(\mff)N}{\mff d}.
\end{equation}
Recall that we have set $\rho=\phi(\mff)/\mff$ in section~\ref{sec3} to simplify typographical work. The reader
will check that the main term (corresponding to $\rho N/[\ell,\ell']$)
vanishes when $q\neq q'$ and is $\rho N/\phi(q)$ otherwise. 
We carry over this change to the bilinear form 
$\bigl\|\sum_q\xi_q\varphi^*_q\bigr\|^2$,
which equals the diagonal term $\rho N\sum_{q}|\xi_q|^2/\phi(q)$
augmented by
\begin{multline*}
  \mathfrak{R}=\sum_{\delta_1\delta_2\delta_3|\Delta}
  \frac{\mu(\delta_2\delta_3)}{\phi(\delta_1)^2\phi(\delta_2)\phi(\delta_3)}
  \sum_{\substack{(\ell,\mathfrak{f}\Delta)=1\\(\ell',\mathfrak{f}\Delta)=1}}
  \frac{\mu(\ell)\xi_{\delta_1\delta_2\ell}}{\phi(\ell)}
  \frac{\mu(\ell')\xi_{\delta_1\delta_3\ell'}}{\phi(\ell')}
  \\\times \sum_{\substack{d|\ell\delta_2\\d'|\ell'\delta_3\\ h|\delta_1}}
  dd'\mu(\ell\delta_2/d)\mu(\ell'\delta_3/d')(\mu\star\phi^2)(h)
  R_{h[d,d']}(M,N,\mathfrak{f}).
\end{multline*}
At this level, we say that
\begin{equation*}
  \bigl|R_{h[d,d']}(M,N,\mathfrak{f})\bigr|\le \theta^*_\mff(h[d,d']/N)
  \le \sum_{0\le r\le R} b_r (h[d,d']/N)^r
\end{equation*}
by~\eqref{upb}. We infer that
\begin{multline*}
  \mathfrak{R}\le
  \sum_{0\le r\le R} b_r N^{-r}
  \sum_{\delta_1\delta_2\delta_3|\Delta}
  \frac{1}{\phi(\delta_1)^2\phi(\delta_2)\phi(\delta_3)}
  \sum_{\substack{(\ell,\mathfrak{f}\Delta)=1\\(\ell',\mathfrak{f}\Delta)=1}}
  \frac{|\xi_{\delta_1\delta_2\ell}|}{\phi(\ell)}
  \frac{|\xi_{\delta_1\delta_3\ell'}|}{\phi(\ell')}
  \\\times \sum_{\substack{d|\ell\delta_2\\d'|\ell'\delta_3\\ h|\delta_1}}
  dd'(\mu\star\phi^2)(h)
  h^r[d,d']^r.
\end{multline*}
This leads to
\begin{multline*}
  \mathfrak{R}\le
  \sum_{0\le r\le R} b_r N^{-r}
  \sum_{\delta_1\delta_2\delta_3|\Delta}
  \frac{\prod_{p|\delta_1}\bigr(1+p^{r+1}(p-2)\bigl)
    \eta_r(\delta_2\delta_3)
  }
  {\phi(\delta_1)^2}
  \\\times 
  \sum_{\substack{(\ell,\mathfrak{f}\Delta)=1\\(\ell',\mathfrak{f}\Delta)=1}}
  |\xi_{\delta_1\delta_2\ell}|\eta_r(\ell)
  |\xi_{\delta_1\delta_3\ell'}|\eta_r(\ell')
  \prod_{\substack{p|(\ell,\ell')}}
  \frac{1+2p^{r+1}+p^{r+2}}{(1+p^{r+1})^2}.
\end{multline*}
The factor that depends on $(\ell,\ell')$ is somewhat troublesome.  We handle
it in the following way: for $r=0$, it is equal to 1; Otherwise let $P$ be the
smallest prime number that does not divide $\mff\Delta$.  This prime
factor will tend to infinity, and we approximate the factor depending on
$(\ell,\ell')$ essentially by $1+\Ocal(P^{-1})$.  More precisely, we write
\begin{multline*}
  \sum_{\substack{(\ell,\mathfrak{f}\Delta)=1\\(\ell',\mathfrak{f}\Delta)=1}}
  |\xi_{\delta_1\delta_2\ell}|\eta_r(\ell)
  |\xi_{\delta_1\delta_3\ell'}|\eta_r(\ell')
  \biggl|\prod_{\substack{p|(\ell,\ell')}}
  \frac{1+2p^{r+1}+p^{r+2}}{(1+p^{r+1})^2}-1\biggr|
  \\\ll_r
  \sum_{p\ge P}
    \sum_{\substack{(m,p\mathfrak{f}\Delta)=1,\\(m',p\mathfrak{f}\Delta)=1}}
    |\xi_{\delta_1\delta_2pm}|\eta_r(pm)
    |\xi_{\delta_1\delta_3pm'}|\eta_r(pm')
  \\\ll_r
    \sum_{p\ge P}p^{2r}
    \sum_{\substack{m,m'}}
    |\xi_{\delta_1\delta_2pm}|\eta_r(m)
    |\xi_{\delta_1\delta_3pm'}|\eta_r(m')
\end{multline*}
The idea here is that the factor $\xi_{\delta_1\delta_2pm}$ forces $m$ to be
rather small. Indeed, anticipating on the values of $\xi$ in \eqref{myval} and
using 
Lemma~\ref{selB}, we get the above to be not more than 
\begin{equation}
  \left(\frac{Z}{\rho N}\right)^2\sum_{p\ge P}p^{2r}
    \left(S/p\right)^{2r+2}\ll_r \left(\frac{Z}{\rho N}\right)^2S^{2r+2}P^{-1}.
\end{equation}
This will give rise to the error term
\begin{equation*}
  \left(\frac{Z}{\rho N}\right)^2
    \sum_{\delta_1\delta_2\delta_3|\Delta}
    \frac{\eta_r^\flat(\delta_1)
    \eta_r(\delta_2\delta_3)}
  {\sigma(\delta_1\delta_2)^{r+1}\sigma(\delta_1\delta_3)^{r+1}}
  \sum_{1\le r\le R} \frac{S^{2r+2}|b_r|}{N^rP}
\end{equation*}
which,  up to a multiplicative constant, is not more than
\begin{equation}
  \left(\frac{Z}{\rho N}\right)^2\prod_{p|\Delta}(1+p^{-1})^2
  \sum_{1\le r\le R} \frac{S^{2r+2}|b_r|}{N^rP}.
\end{equation}
The factor $P^{-1}$ will indeed be enough to control this quantity.
Hence, again anticipating  on~\eqref{myval}, we reach
\begin{multline*}
  \Bigl\|\sum_q\xi_q\varphi^*_q\Bigr\|^2
  \le \rho N\sum_{q}|\xi_q|^2/\phi(q)
  \\+\sum_{0\le r\le R} \frac{b_r}{N^r}
  \sum_{\delta_1\delta_2\delta_3|\Delta}
  \eta_r^\flat(\delta_1)
  \sum_{\substack{(\ell,\mathfrak{f}\Delta)=1,\\(\ell',\mathfrak{f}\Delta)=1}}
    |\xi_{\delta_1\delta_2\ell}|\eta_r(\delta_2\ell)
    |\xi_{\delta_1\delta_3\ell'}|\eta_r(\delta_3\ell')
  \\+\Ocal\biggl( \left(\frac{Z}{\rho N}\right)^2\prod_{p|\Delta}(1+p^{-1})^2
  \sum_{1\le r\le R} \frac{S^{2r+2}|b_r|}{N^rP}
    \biggr).
\end{multline*}

\section{Some arithmetical auxiliaries}

We need to evaluate some rather unusual averages. 
\begin{lem}\label{selB} Let $\mff^*$ be a positive integer. We set
  $\rho^*=\phi(\mff^*)/\mff^*$ and use $t(q)=1-\sigma(q)/S^*$. For
  any real number $S^*$ going to infinity, we have
  \begin{equation*}
    \sum_{\substack{q/\sigma(q)\le S^*,\\(q,\mff^*)=1}}\frac{t(q)^2}{\phi(q)}
    =\rho^*(\Log S^*+\kappa(\mff^*))+o(1)
  \end{equation*}
  with
  \begin{equation*}
    \kappa(\mff)=\gamma+\sum_{p\ge2}\frac{\Log p}{p(p-1)}
    -\sum_{p}\frac{\Log(1+p^{-1})}{p}
    +\sum_{p|\mathfrak{f}^*}\frac{\Log (p+1)}{p}
    -\tfrac32
  \end{equation*}
  ($\kappa(210)=1.115\,37\dots$)
  and
  \begin{equation*}
    \sum_{\substack{q/\sigma(q)\le S^*,\\(q,\mff^*)=1}}\eta_r(q)t(q)
    =
    \frac{\rho^*}{2(r+1)}
    \prod_{p\nmid\mff^*}\left(1-\frac{p^r-1}{p^{r+1}(p+1)}\right)
    S^{*(r+1)}+o(S^{*(r+1)}).
  \end{equation*}
\end{lem}

\begin{dem}
  The first estimate comes from~\cite{Selberg*91}.
  We follow closely Selberg's proof and get
  \begin{equation*}
    \sum_{\substack{q/\sigma(q)\le
        S^*,\\(q,\mff^*)=1}}\frac{\eta_r(q)t(q)}{q^r}
    =\frac{\rho}{2}\prod_{p\nmid\mff^*}\left(1-\frac{p^r-1}{p^{r+1}(p+1)}\right)
    S^*+o(S^*)
  \end{equation*}
  from which we get
  \begin{equation*}
    \sum_{\substack{q/\sigma(q)\le
        S^*,\\(q,\mff^*)=1}}\eta_r(q)t(q)
    =\frac{\rho^*}{2(r+1)}
    \prod_{p\nmid\mff^*}\left(1-\frac{p^r-1}{p^{r+1}(p+1)}\right)
    S^{*(r+1)}+o(S^{*(r+1)}).
  \end{equation*}
\end{dem}

Note that the quantities we end up computing are the same as the ones that appear
in~\cite{Selberg*91} though we have one less
to handle.

Let us define
  \begin{equation}
    C_r(\Delta)=
    \frac{\phi(\Delta)^2}{\Delta^2}
    \sum_{\delta_1\delta_2\delta_3|\Delta}
  \frac{\prod_{p|\delta_1}\bigr(1+p^{r+1}(p-2)\bigl)
    \eta_r(\delta_2\delta_3)}
  {\phi(\delta_1)^2
    \sigma(\delta_1)^{2r+2}\sigma(\delta_2)^{r+1}\sigma(\delta_3)^{r+1}}.
\end{equation}

We have
\begin{lem}
\begin{multline*}
C_r(\Delta) = \prod_{p|\Delta} \biggl(\frac{(p-1)^2}{p^2} +
\frac{2(p-1)(p^{r+1}+1)}{p^2(p+1)^{r+1}}\\
+\frac{1+p^{r+1}(p-2)}{p^2(p+1)^{2r+2}} +
\frac{(1+p^{r+1}(p-2))(p^{r+1}+1)}{(p-1)p^2(p+1)^{3r+3})}\biggr).
\end{multline*}
\end{lem}


\begin{dem}
  We start with $\delta_3$:
  \begin{eqnarray*}
    \sum_{\delta_3|\Delta/(\delta_1\delta_2)}
  \frac{
    \eta_r(\delta_3)}
  {
    \sigma(\delta_3)^{r+1}}
  &=&\prod_{p|\Delta/(\delta_1\delta_2)}
  \left(1+\frac{1+p^{r+1}}{(p+1)^{r+1}(p-1)}\right)
  \\&=&\prod_{p|\Delta/(\delta_1\delta_2)}
  \frac{(p+1)^{r+1}(p-1)+p^{r+1}+1}{(p+1)^{r+1}(p-1)}.
  \end{eqnarray*}
  Our sum reduces to
  \begin{multline*}
    \prod_{p|\Delta}
  \frac{\Bigl((p+1)^{r+1}(p-1)+p^{r+1}+1\Bigr)(p-1)}{p^{2}(p+1)^{r+1}}
  \\\times\sum_{\delta_1\delta_2|\Delta}
  \frac{\prod_{p|\delta_1}\bigr(1+p^{r+1}(p-2)\bigl)
    \eta_r(\delta_2)}
  {\phi(\delta_1)^2
    \sigma(\delta_1)^{2r+2}\sigma(\delta_2)^{r+1}}\prod_{p|\delta_1\delta_2}
  \frac{(p+1)^{r+1}(p-1)}{(p+1)^{r+1}(p-1)+p^{r+1}+1}.
  \end{multline*}
  We continue with $\delta_2$:
  \begin{eqnarray*}
    \sum_{\delta_2|\Delta/\delta_1}
  &&\frac{1+p^{r+1}}{(p-1)(p+1)^{r+1}}
  \frac{(p+1)^{r+1}(p-1)}{(p+1)^{r+1}(p-1)+p^{r+1}+1}
  \\&=&\prod_{p|\Delta/\delta_1}
  \left(
    1
    +
    \frac{p^{r+1}+1}{(p+1)^{r+1}(p-1)+p^{r+1}+1}
  \right)
  \\&=&\prod_{p|\Delta/\delta_1}
    \frac{(p+1)^{r+1}(p-1)+2p^{r+1}+2}{(p+1)^{r+1}(p-1)+p^{r+1}+1}.
  \end{eqnarray*}
  Hence our quantity reduces to
  \begin{multline*}
    \prod_{p|\Delta}
  \frac{\bigl((p+1)^{r+1}(p-1)+2p^{r+1}+2\bigr)(p-1)}{p^2(p+1)^{r+1}}
  \\\times\sum_{\delta_1|\Delta}\prod_{p|\delta_1}\frac{1+p^{r+1}(p-2)}
  {(p-1)^2(p+1)^{2r+2}}
  \frac{(p+1)^{r+1}(p-1)+p^{r+1}+1}{(p+1)^{r+1}(p-1)+2p^{r+1}+2}
\end{multline*}
which reads
  \begin{multline*}
  \prod_{p|\Delta}
  \frac{\bigl((p+1)^{r+1}(p-1)+2p^{r+1}+2\bigr)(p-1)}{p^2(p+1)^{r+1}}
  \\\times\frac{\begin{array}{r}(p-1)^2(p+1)^{2r+2}\bigl((p+1)^{r+1}(p-1)+2p^{r+1}+2\bigr)
      \\ \ \qquad+\bigl(1+p^{r+1}(p-2)\bigr)\bigl((p+1)^{r+1}(p-1)+p^{r+1}+1\bigr)
    \end{array}}{(p-1)^2(p+1)^{2r+2}\bigl((p+1)^{r+1}(p-1)+2p^{r+1}+2\bigr)}
  \\=
  \prod_{p|\Delta}
  \frac{\begin{array}{c}(p-1)^2(p+1)^{2r+2}\bigl((p+1)^{r+1}(p-1)+2p^{r+1}+2\bigr)
  \\ \ +\bigl(1+p^{r+1}(p-2)\bigr)\bigl((p+1)^{r+1}(p-1)+p^{r+1}+1\bigr)\end{array}}{(p-1)p^2(p+1)^{3r+3}}.
\end{multline*}
\end{dem}

\section{Using the hermitian inequalilty}

Optimizing in $\xi$ is too difficult. We stick to the simplest choice:
$M_i=\rho N/\phi(q)$, $[f|\varphi_i^*]/M_i=Z/(\rho N)$,
$n_i=\sigma(q)/\phi(q)$ 
and $Y=Z/S$ for a parameter $S$ we shall choose later on. This leads to
\begin{equation}\label{myval}
  \xi_q=\frac{Z}{\rho N}t(q),\quad t(q)=1-\frac{\sigma(q)}{S}.
\end{equation}

We invoke Lemma~\ref{selB} to compute the relevant mean values, and for
instance, we use $S^*=S/\sigma(\delta_1\delta_2)$ and
$\mff^*=\mathfrak{f}\Delta$ to evaluate $\sum_{\substack{(\ell,\mathfrak{f}\Delta)=1}}
  |\xi_{\delta_1\delta_2\ell}|\eta_r(\ell)$.
  There
  appear  constants in the form of an Euler product, say
  $\mathfrak{S}_r(\mff^*)$, which we again approximate 
  by $1+\Ocal(P^{-1})$. Let us give some details.
In a first step we reach
\begin{equation*}
  Z\ge\left(\frac{Z}{\rho N}\right)^2\rho^2N
  \left(\Log S+\kappa(\mff)\right)
  +\frac{2Z^2}{\rho N S}\sum_{(q,\mff)=1} \frac{\sigma(q)t(q)}{\phi(q)}
  -A+\Ocal(B)
\end{equation*}
with $g(\delta)=\prod_{p|\delta}(1+p^{r+1}(p-2))/(p-1)^2$ and 
\begin{equation*}
  \left\{
    \begin{array}{l}
      \displaystyle 
      A=\sum_{0\le r\le R} \frac{Z^2b_r}{\rho^2N^{r+2}}
  \sum_{\delta_1\delta_2\delta_3|\Delta}
  g(\delta_1)
  \sum_{\substack{(\ell,\mathfrak{f}\Delta)=1\\(\ell',\mathfrak{f}\Delta)=1}}
    t(\delta_1\delta_2\ell)\eta_r(\ell\delta_2)
    t(\delta_1\delta_3\ell')\eta_r(\ell'\delta_3),
    \\ \displaystyle 
    B=\sum_{\delta_1\delta_2\delta_3|\Delta}
  g(\delta_1)
    \sum_{1\le r\le R} \sum_{p\ge P}\frac {|b_r|p^{2r}}{N^r}
    \sum_{\substack{\ell,\ell'}}
    |\xi_{\delta_1\delta_2p\ell}|\eta_r(\ell\delta_2)
    |\xi_{\delta_1\delta_3p\ell'}|\eta_r(\ell'\delta_3).
    \end{array}
  \right.
\end{equation*}
We tidy this expression step by step:
\begin{multline*}
  N\ge Z
  \left(\Log S+\kappa(\mff)\right)
  +Z 
    \\-Z\sum_{0\le r\le R} \frac{b_r  S^{2r+2}}{\rho^2N^{r+1}}
  \sum_{\delta_1\delta_2\delta_3|\Delta}
  \frac{\prod_{p|\delta_1}\bigr(1+p^{r+1}(p-2)\bigl)
    \eta_r(\delta_2\delta_3)}
  {\phi(\delta_1)^2
    \sigma(\delta_1)^{2r+2}\sigma(\delta_2)^{r+1}\sigma(\delta_3)^{r+1}}
  \frac{\rho^2\phi(\Delta)^2\mathfrak{S}_r(\mff\Delta)^2}{4\Delta^2(r+1)^2}
  \\+\Ocal\Bigl(\prod_{p|\Delta}(1+p^{-1})^2S^2ZN^{-1}P^{-1}
    \sum_{1\le r\le R} |b_r|(S^2/N)^r
    \Bigr)+o(N).
\end{multline*}
This leads to
\begin{multline*}
  N/Z\ge 
  \Log S+\kappa(\mff)
  +1
    -
    \sum_{0\le r\le R} \frac{b_r (S^2/N)^{r+1}}{4(r+1)^2}
  C_r(\Delta)
  \mathfrak{S}_r(\mff\Delta)^2
  \\+\Ocal\Bigl(\prod_{p|\Delta}(1+p^{-1})^2P^{-1}
    \sum_{1\le r\le R} |b_r|(S^2/N)^{r+1}
    \Bigr)+o(1)
\end{multline*}
And since $\mathfrak{S}_r(\mff\Delta)=1+\Ocal(P^{-1})$, we finally obtain
\begin{multline*}
  N/Z-\tfrac12\Log N\ge 
  \tfrac12\Log (S^2/N)+\kappa(\mff)
  +1
    -
    \sum_{0\le r\le R} \frac{b_r (S^2/N)^{r+1}}{4(r+1)^2}
  C_r(\infty/\mff)
  \\+\Ocal\Bigl(\prod_{p|\Delta}(1+p^{-1})^2P^{-1}
    \sum_{1\le r\le R} |b_r|(S^2/N)^{r+1}
    \Bigr)+o(1).
\end{multline*}
At this level, we send $\Delta$ (and $P$) to infinity and we are left with
finding an optimal value for $S^2/N$. It would be satisfactory to have an
expression for the final constant, but we are not able to attain such
a precision. In particular, the $b_r$'s should not appear in the final
expression. We are, however, able to get numerical results.




\section{Optimizing the polynomial via linear programming}
\label{otp}

It is better at this level to change notation slightly and set
\begin{multline*}
c_r(x) = \prod_{7< p\leq x} \biggl(\frac{(p-1)^2}{p^2} +
\frac{2(p-1)(p^{r+1}+1)}{p^2(p+1)^{r+1}}\\
+\frac{1+p^{r+1}(p-2)}{p^2(p+1)^{2r+2}} +
\frac{(1+p^{r+1}(p-2))(p^{r+1}+1)}{(p-1)p^2(p+1)^{3r+3})}\biggr),
\end{multline*}
together with
\begin{equation}
  \label{eq:3}
  c_r =C_r(\infty/\mff) = \lim_{x\to\infty} c_r(x).
\end{equation}
\begin{lem}
We have for $x>2r+2$ the estimate
\[
c_r(x)\geq c_r \geq c_r(x)\Bigl(1-\frac{2r+2}{x}\Bigr).
\]
\end{lem}
\begin{dem}
Denote by $F(p)$ the factor in $c_r(x)$ corresponding to the prime
number $p$. We then have
\begin{multline*}
\frac{(p-1)^2}{p^2} +
\frac{2(p-1)(p^{r+1}+1)}{p^2(p+1)^{r+1}} \leq F(p) \\\leq \frac{(p-1)^2}{p^2} +
\frac{2(p-1)(p^{r+1}+1)}{p^2(p+1)^{r+1}} + \frac{2}{p^2(p+1)^r},
\end{multline*}
and the left-hand expression can be written as
\[
1 - \frac{2 - 2p + 2p^{1+r} - 2p^{2+r} - (1+p)^r + 
p(1+p)^r + 2p^2(1+p)^r}{p^2(1+p)^{r+1}},
\]
from which it is obvious that $F(p)<1$ for $r\geq 1$. This can
be checked directly when $r=0$. On the other hand we have for $p\geq r$ the estimate
\[
p^r =((p+1)-1)^r = \sum_{\nu=0}^r (-1)^\nu\binom{r}{\nu} (p+1)^{r-\nu} \geq (p+1)^r
- r(p+1)^{r-1},
\]
since the binomial sum is alternating and monotonically decreasing,
which implies, for $r\geq 1$, the estimate
\begin{eqnarray*}
F(p) & \geq & 1 - \frac{2 - 2p + 2p^{1+r} - (1+p)^r + 
p(1+p)^r + 2rp^2(1+p)^{r-1}}{p^2(1+p)^{r+1}}\\
 & \geq & 1-\frac{(2r+2)p^2(p+1)^{r-1}}{p^2(p+1)^{r+1}}\\
 & \geq & \left(1-\frac{1}{p^2}\right)^{2r+2};
\end{eqnarray*}
Again, for $r=0$ this bound can be checked directly. We now have
\begin{multline*}
c_r=c_r(x)\prod_{p>x} F(p) \geq c_r(x) \prod_{p>x}
\left(1-\frac{1}{p^2}\right)^{2r+2} \geq c_r(x)\Bigl(1-\sum_{n>x}
  n^{-2}\Bigr)^{2r+2}\\
 \geq c_r(x) \Bigl(1-\frac{1}{x}\Bigr)^{2r+2}
 \geq c_r(x)\Bigl(1-\frac{2r+2}{x}\Bigr).
\end{multline*}
\end{dem}

Next we compute a polynomial $P$ such that $P(V n/300)>\theta^*_{210}(V n/300)$ for
all $n\leq 300$, where $V=S^2/N$ is some parameter, and such that the linear
functional $\mathcal{F}^*(P)$ defined on monomials as
$\mathcal{F}^*(x^r)=\frac{1}{4(r+1)^2}c_r(10^5)$ is minimized. Since the
domain described by the inequalities relating $P$ and $\theta^*_{210}$ is
non-compact, we further restrict all coefficients to be at
least $-M$. Then we compute an upper bound for $\mathcal{F}(P)$, the
linear functional defined on monomials through
$\mathcal{F}(x^r)=\frac{1}{4(r+1)^2}c_r$, by replacing $c_r$ by
$c_r(1-(2r+2)/10^5)$, whenever
the coefficient of $x^r$ in $P$ is negative. Note that smaller values
of $M$ yield worse approximations to the optimal value of the
functional $\mathcal{F}^*(P)$, but restricting the size of negative
coefficients diminishes the upper bound for $\mathcal{F}(P)$. 

Using $M=-1\,000$, $V=1.5$ and looking for polynomials of degree~25, we
obtain a polynomial $P$ with $\mathcal{F}(P)\leq 0.547\,38$.
The resulting polynomial is 
\begin{multline*}
P(u) = 3.117973+ 3.555433\,u-154.037413 \,u^{2} + 732.936467 \,u^{3} -1000
\,u^{4} 
\\-1000 \,u^{5} +
3227.305717 \,u^{6} -1000 \,u^{7} -1000 \,u^{8} -1000 \,u^{9} -1000 \,u^{10} 
\\+ 3012.745710 \,u^{11} +
1227.721539 \,u^{12} -1000 \,u^{13} -1000 \,u^{14} 
-1000 \,u^{15} \\-1000 \,u^{16} + 1191.986883 \,u^{17} +
2708.564854 \,u^{18} -1000 \,u^{19} 
-1000 \,u^{20} 
\\-1000 \,u^{21} + 675.282733 \,u^{22} +
1158.017142 \,u^{23} -1000 \,u^{24} + 214.336183 u^{25};
\end{multline*}

We determined this polynomial by using the lpsolve linear programming package
and a C-program of our own. There has been numerous precisions issues and
instabilities that we were no able to understand, less to tackle to our satisfaction.
For instance, many of the coefficients are on the
artificial boundary $a_i\geq -1\,000$. We delay a further study to a latter
paper.

Once the polynomial is selected, we can simply consider it and study anew how
it fits $\theta^*_\mff$. To do so, we revert to Pari/GP.

By construction of $P$ we know that $P(u)\geq\theta^*_\mff(u)$ for 300 well spaced
points, however, this does not imply that this inequality holds true
for all values of $u$. In fact, there are six regions in which $P$ dips
slightly below $\theta^*_{210}$, these regions being close to the points $u=0.29,
0.59, 0.84, 1.16, 1.32$ and $1.44$. The difference is greatest at $u=1.442\,618\dots$,
where $P(u)-\theta^*_{210}(u) = -0.008\,338\dots$. Hence, putting $P^*(u)=P(u)+0.0084$,
we obtain a polynomial which is strictly larger than $\theta^*_{210}$, and we have
\[
\mathcal{F}(P^*) = \mathcal{F}(P) + \frac{0.0084}{4} \leq 0.547\,4+
0.002\,1 = 0.549\,5,
\]
which implies that
\begin{equation*}
N/Z-\frac{1}{2}\Log N  \geq  \tfrac{1}{2}\Log 1.5 + \kappa(210) + 1 -
\mathcal{F}(P^*)
  \geq  1.768\,6
\end{equation*}
and therefore
\[
\pi(M+N) - \pi(M) \leq \frac{2N}{\Log N + 3.537\,2}
\]
as announced.


\begin{minipage}[t]{0.5\linewidth}
  Olivier Ramar\'e

  Laboratoire Paul Painlev\'e
 
  Universit\'e Lille 1

  59 655 Villeneuve d'Ascq, cedex

  France

  \texttt{ramare@math.univ-lille1.fr}
\end{minipage}
\begin{minipage}[t]{0.5\linewidth}
  Jan-Christoph Schlage-Puchta

  Mathematisches Institut Freiburg

  Freiburg

  Germany
\end{minipage}

\end{document}